\documentclass[]{article}

\title{On finite groups with all simple modules of low dimension in characteristic $p$}
\author{Geoffrey R. Robinson}

\usepackage{latexsym,amssymb}
\begin{document}

\maketitle

\begin{abstract}
	We offer a short and reasonably elementary proof that if $G$ is a finite group, $F$ is an algebraically closed field of prime characteristic $p$, and all simple $FG$-modules have dimension less than $p$, then $G$ has a normal Sylow $p$-subgroup. The algebraic closure of $F$ is not really necessary, but simplifies the exposition in the essential case.

\end{abstract}

\section{Introduction}

\medskip
Let $G$ be a finite group, $p$ be a prime, and let $F$ be an algebraically closed field of characteristic $p$. We recall that $O_{p}(G)$ acts trivially on each simple $FG$-module. We will prove that if each simple $FG$-module has dimension less than $p$, then $G$ has a normal Sylow $p$-subgroup. This may be compared with an analogous result for ordinary irreducible representations proved in [2] by Isaacs and Passman, where the resulting normal Sylow $p$-subgroup is necessarily Abelian, which need not be the case in the present modular context.

\medskip
We will make use of the following reciprocity theorem, which appears in [3] with a proof using properties of the Reynolds ideal of $Z(FG)$, but which may also be proved using projective homomorphisms (see, for example, Lemma 5.2 of Grodal [1]): if $H$ is a subgroup of $G$ and $S$ is simple $FG$-module, $T$ is a simple $FH$-module, then the multiplicity of the projective cover of $S$ as a summand of ${\rm Ind}^{G}_{H}(T)$ is equal to the multiplicity of the projective cover of $T$ as a summand of ${\rm Res}^{G}_{H}(S).$ 

\medskip
Now we prove our main: 

\medskip
\noindent {\bf Theorem:} \emph{ Let $G$ be a finite group and $F$ be an algebraically closed field of prime characteristic $p$. Suppose that every simple $FG$-module has dimension less than $p$.
	Then $G$ has a normal Sylow $p$-subgroup.}

\medskip
\noindent {\bf Proof :} We proceed by induction on $|G|$. Since $O_{p}(G)$ acts trivially on every simple $FG$-module, we may suppose that $O_{p}(G) = 1$, but that $p$ divides $|G|$.

\medskip
Since the hypotheses imply (on consideration of the composition factors of the respective regular modules) that for any section $K$ of $G$, each simple $FK$-module has dimension less than $p$, we may suppose by induction that every proper section $K$ of $G$ has a normal Sylow $p$-subgroup. Let $P$ be a Sylow $p$-subgroup of $G$.

\medskip
We claim that $P \cap P^{g} = 1$ for all $g \in G \backslash N_{G}(P)$. It suffices to prove that $N_{G}(V) \leq N_{G}(P)$ whenever $1 \neq V \leq P$. Since $O_{p}(G) = 1$, we have 
$N_{G}(V) < G$ for all such $V$. Hence, by induction, we know that $N_{G}(V)$ has a normal (and hence unique) Sylow $p$-subgroup for each such $V$. It follows that
$N_{G}(V)$ is contained in the normalizer of its unique Sylow $p$-subgroup, so, using induction on $[P:V]$ (and the fact that $N_{P}(V) > V$ whenever $V < P$), the claim follows.

\medskip
Now  $N_{G}(P) < G$ and we have ${\rm Ind}_{N_{G}(P)}^{G}(F) = F \oplus  M$ where $M$ is a (non-zero) projective $FG$-module, using Mackey decomposition, for example.  Let $Q$ be a projective indecomposable summand of $M$, and let $S$ be the socle of $Q$ (which is also isomorphic to the head of $Q$). Then the multiplicity of the projective cover of the trivial $FN_{G}(P)$-module $F$ as a summand of  ${\rm Res}^{G}_{N_{G}(P)}(S)$ is strictly positive by the above reciprocity theorem (applied with the trivial $FN_{G}(P)$-module $F$ in place of $T$). Hence ${\rm dim}_{F}(S) \geq |P|$, contrary to our assumption that  each simple $FG$-module has dimension less than $p$. This contradiction shows that $|G|$ is not divisible by $p$, and so completes the proof of the Theorem.

\medskip

\begin{center}
	{\bf References}
\end{center}	

\medskip
\noindent [1] Grodal, J. ; \emph{Higher limits via subgroup complexes}, Annals of Math. (2), {\bf 155} (2),2002,405-457.

\medskip
\noindent [2] Isaacs, I.M and Passman, D.S. ; \emph{A characterization of groups in terms of the degree of their characters}, Pacific J. Math., {\bf 15}, 3, (1965),877-903.

\medskip
\noindent [3] Robinson, G.R.; \emph{On projective summands of induced modules}, Journal of Algebra, {\bf 122}, (1989),106-111.

\end{document}